\documentclass{amsart}
\usepackage{upref}
\usepackage{amscd}
\numberwithin{equation}{section}
\newtheorem{theorem}[equation]{Theorem}
\newtheorem{lemma}[equation]{Lemma}
\newtheorem{corollary}[equation]{Corollary}
\theoremstyle{definition}
\newtheorem{remark}[equation]{Remark}
\newtheorem{definition}[equation]{Definition}
\DeclareMathOperator{\im}{im}
\newcommand{\Smash}{\wedge}
\newcommand{\dual}[1]{D \mspace{-2mu} {#1}}
\newcommand{\invlim}{\varprojlim}
\newcommand{\conn}[1]{\lvert #1 \rvert}
\newcommand{\mathcolon}{\colon\,}
\newcommand{\exactcouple}[6]
 {\begin{picture}(16,8)(0,1.7)	   %% maps go counterclockwise.  args:
  \put(2,8){\ensuremath{#1}}	   %% top left, bottom, top right,
  \put(8,2){\ensuremath{#2}}	   %% bottom left map, bottom right
  \put(14,8){\ensuremath{#3}}	   %% map, top map
  \put(5.2,4.8){\makebox(0,0)[r]{\ensuremath{#4}}}
  \put(11.8,4.8){\makebox(0,0)[l]{\ensuremath{#5}}}
  \put(8.5,9){\makebox(0,0){\ensuremath{#6}}}
  \put(3.5,7.5){\vector(1,-1){4.5}}
  \put(9,3){\vector(1,1){4.5}}
  \put(13.5,8.2){\vector(-1,0){9.5}}
 \end{picture}}
\begin{document}

\title[Vanishing lines]
{Vanishing lines in generalized Adams spectral sequences are generic}

\date{\today}

\author{M. J. Hopkins}
\address{Department of Mathematics \\ 
Massachusetts Institute of Technology \\
Cambridge, MA 02139}
\email{mjh@math.mit.edu}
\author{J. H. Palmieri}
\address{Department of Mathematics \\ 
University of Notre Dame \\
Notre Dame, IN 46556}
\email{palmieri@member.ams.org}
\author{J. H. Smith}
\address{Department of Mathematics \\ 
Purdue University \\
West Lafayette, IN 47907}
\email{jhs@math.purdue.edu}


 % stable homotopy theory, spectra

\begin{abstract}
We show that in the generalized Adams spectral sequence, the presence
of a vanishing line of fixed slope at some $E_{r}$-term is a generic
property.
\end{abstract}

\maketitle

\section{Introduction}
\label{sec-intro}

Let $E$ be a nice ring spectrum, let $X$ be a spectrum, and consider
the $E$-based Adams spectral sequence converging to $\pi_{*}X$.  In
this note, we prove that, for any number $m$, the property that the
spectral sequence has a vanishing line of slope $m$ at some term of
the spectral sequence is generic.

\begin{definition}
We say that a spectrum $X$ is \emph{$E$-complete} if the inverse limit
of the Adams tower for $X$ is contractible.  (See
Section~\ref{sec-proof} for a definition of the Adams tower.)

A property $P$ of $E$-complete spectra is said to be \emph{generic} if
\begin{itemize}
\item whenever $Y$ is $E$-complete and $Y$ satisfies $P$, then so does
any retract of $Y$; and
\item if $X \rightarrow Y \rightarrow Z$ is a cofibration of
$E$-complete spectra and two of $X$, $Y$, and $Z$ satisfy $P$, then so
does the third.
\end{itemize}
In other words, a property is generic if the full subcategory of all
$E$-complete spectra satisfying it is thick.
\end{definition}

Given a connective spectrum $W$, we write $\conn{W}$ for its
connectivity.

We assume that our ring spectrum $E$ satisfies the standard
assumptions for convergence of the $E$-based Adams spectral
sequence---in other words, the assumptions necessary for Theorem
15.1(iii) in \cite[Part III]{adams-blue}; see also
Assumptions~2.2.5(a)--(c) and (e) in \cite{ravenel-green}.

\begin{theorem}\label{thm-main}
Let $E$ be a ring spectrum as above, and consider the $E$-based Adams
spectral sequence $E_{*}^{**}(X) \Rightarrow \pi_{*}(X)$.  Fix a
number $m$.  The following properties of an $E$-complete spectrum $X$
are each generic\textup{:}
\begin{itemize}
\item [(i)] There exist numbers $r$ and $b$ so that for all $s$ and $t$
with $s \geq m(t-s) + b$, we have $E_{r}^{s,t}(X) = 0$.
\item [(ii)] There exist numbers $r$ and $b$ so that for all finite
spectra $W$ with $\conn{W} = w$ and for all $s$ and $t$ with $s \geq
m(t-s-w)+b$, we have $E_{r}^{s,t}(X \Smash W) = 0$.
\end{itemize}
\end{theorem}

\begin{remark}\label{remark-main}
\begin{itemize}
\item [(a)] One usually draws Adams spectral sequences $E_{r}^{s,t}$
with $s$ on the vertical axis and $t-s$ on the horizontal; in terms of
these coordinates, the properties say that $E_{r}^{s,t}$ is zero above
a line of slope $m$.
\item [(b)] Assuming that $X$ is $E$-complete ensures that the
spectral sequence converges, which we need to prove the theorem.  We
do not need to identify the $E_{2}$-term of the spectral sequence, so
we do not need to know that $E$ is a flat ring spectrum, for example.
\end{itemize}
\end{remark}

We also mention one or two possible applications of the theorem.
Since there is a classification of the thick subcategories of the
category of finite spectra (see \cite{hopkins-glob, hopkins-smith,
ravenel-orange}), then if one is dealing with finite spectra $X$, one
may be able to identify all spectra with vanishing line of a given
slope.  For example, in the classical mod $2$ Adams spectral sequence,
since the mod $2$ Moore spectrum has a vanishing line of slope
$\frac{1}{2}$ at the $E_{2}$-term, then the mod $2^{n}$ Moore
spectrum, and indeed any type 1 spectrum, has a vanishing line of
slope $\frac{1}{2}$ at some $E_{r}$-term.  Similarly, any type $n$
spectrum has a vanishing line of slope $\frac{1}{2^{n+1}-2}$ at some
$E_{r}$-term of the classical mod 2 Adams spectral sequence.
Theorem~\ref{thm-main} gives no control over the term $r$ or the
intercept $b$ of the vanishing line.

Since the proof is formal, this theorem also applies in other stable
homotopy categories.  The second author has used this result in an
appropriate category of modules over the Steenrod algebra to prove a
version of Quillen stratification for the cohomology of the Steenrod
algebra.  See \cite{palmieri-f-iso, palmieri-steenrod} for
details.

\section{Proof of Theorem~\ref{thm-main}}
\label{sec-proof}

The difficulty in proving a result like Theorem~\ref{thm-main} is that
the $E_{r}$-term of an Adams spectral sequence does not have nice
exactness properties if $r \geq 3$---a cofibration of spectra does not
lead to a long exact sequence of $E_{r}$-terms, for instance.  So we
prove the theorem by showing that the purported generic conditions are
equivalent to other conditions on composites of maps in the Adams
tower, and then we show that those other conditions are generic.

We start by describing the standard construction of the Adams spectral
sequence, as found in \cite[III.15]{adams-blue},
\cite[2.2]{ravenel-green}, and any number of other places.  Given a
ring spectrum $E$, we let $\overline{E}$ denote the fiber of the unit
map $S^{0} \xrightarrow{} E$.  For any integer $s \geq 0$, we let
\begin{gather*}
F_{s} X = \overline{E}^{\Smash s} \Smash X, \\
K_{s} X = E \Smash \overline{E}^{\Smash s} \Smash X.
\end{gather*}
We use these to construct the following diagram of cofibrations, which
we call the \emph{Adams tower for $X$}:
\[
\begin{CD}
X @= F_{0} X @<{g}<< F_{1} X @<{g}<< F_{2} X @<{g}<< \dots. \\
@. @VVV @VVV @VVV \\
@. K_{0} X @. K_{1} X @. K_{2} X
\end{CD}
\]
This construction satisfies the definition of an ``$E_{*}$-Adams
resolution'' for $X$, as given in \cite[2.2.1]{ravenel-green}---see
\cite[2.2.9]{ravenel-green}.  Note also that $F_{s} X = X \Smash F_{s}
S^{0}$, and the same holds for $K_{s} X$---the Adams tower is
functorial and exact.

Given the Adams tower for $X$, if we apply $\pi_{*}$, we get an exact
couple and hence a spectral sequence.  This is called the
\emph{$E$-based Adams spectral sequence}.  More precisely, we let
\begin{gather*}
D_{1}^{s,t} = \pi_{t-s} F_{s} X, \\
E_{1}^{s,t} = \pi_{t-s} K_{s} X.
\end{gather*}
If we let $g \mathcolon F_{s+1}X \xrightarrow{} F_{s}X$ denote the
natural map, then $g_{*} = \pi_{t-s}(g)$ is the map $D_{1}^{s+1,t+1}
\xrightarrow{} D_{1}^{s,t}$.  Then we have the following exact couple
(the pairs of numbers indicate the bidegrees of the maps):
\begin{center}
\exactcouple{D_{1}^{s,t}}{E_{1}^{s,t}}{D_{1}^{s+1,t+1}}
   {(0,0)}{(1,0)}{(-1,-1)}
\end{center}
This leads to the following $r$th derived exact couple, where
$D_{r}^{s,t}$ is the image of $g_{*}^{r-1}$, and the map
$D_{r}^{s+1,t+1} \xrightarrow{} D_{r}^{s,t}$ is the restriction of
$g_{*}$:
\begin{center}
\exactcouple{D_{r}^{s,t}}{E_{r}^{s+r-1,t+r-1}}
 {D_{r}^{s+1,t+1}} 
 {(r-1,r-1)}{(1,0)}{(-1,-1)} 
\end{center}
Unfolding this exact couple leads to the following exact
sequence:
\begin{equation}\label{eqn-ass-exact}
\dots \xrightarrow{} E_{r}^{s,t+1} \xrightarrow{} D_{r}^{s+1,t+1}
\xrightarrow{} D_{r}^{s,t} \xrightarrow{} E_{r}^{s+r-1,t+r-1}
\xrightarrow{} \dots .
\end{equation}

Fix a number $m$.  With respect to the $E$-based Adams spectral
sequence $E_{*}^{**}(-)$, we have the following conditions on a
spectrum $X$:
\begin{itemize}
\item [(1)] There exist numbers $r$ and $b$ so that for all $s$ and
$t$ with $s \geq m(t-s) + b$, the map $ g_{*}^{r-1} \mathcolon
\pi_{t-s}(F_{s+r-1}X) \xrightarrow{} \pi_{t-s}(F_{s}X)$ is zero.  (In
other words, $D_{r}^{s,t}(X) = 0$.)
\item [(2)] There exist numbers $r$ and $b$ so that for all $s$ and
$t$ with $s \geq m(t-s) + b$, we have $E_{r}^{s,t}(X) = 0$.
\item [(3)] There exist numbers $r$ and $b$ so that for all finite
spectra $W$ with $\conn{\dual{W}} = -w$ and for all $s$ with $s \geq
mw + b$, then the composite $W \xrightarrow{} F_{s+r-1}X
\xrightarrow{} F_{s}X$ is null.  (Here, $\dual{W}$ denotes the
Spanier-Whitehead dual of $W$.)
\item [(4)] There exist numbers $r$ and $b$ so that for all finite
spectra $W$ with $\conn{W} = w$ and for all $s$ and $t$ with $s \geq
m(t-s-w) + b$, we have $E_{r}^{s,t}(X \Smash W) = 0$.
\end{itemize}
Each condition depends on a pair of numbers $r$ and $b$, and we write
$(1)_{r,b}$ to mean that condition (1) holds with the numbers
specified, and so forth.

Notice that if $m=0$, then condition (3) says that $F_{s+r-1}X
\xrightarrow{} F_{s}X$ is a phantom map whenever $s \geq b$.  If
$m=0$, then condition (1) says that $F_{s+r-1}X \xrightarrow{} F_{s}X$
is a ghost map (zero on homotopy) whenever $s \geq b$.

\begin{lemma}\label{lemma-main}
Fix numbers $m$, $r$, and $b$.  We have the following
implications:
\begin{itemize}
\item [(a)] If $r \geq -m$, then $(1)_{r,b} \Rightarrow
(2)_{r,b+r-1}$.  If $r < -m$, then $(1)_{r,b} \Rightarrow
(2)_{r,b-m}$.
\item [(b)] If $r \geq 1-m$, then $(2)_{r,b} \Rightarrow
(1)_{r,b-m}$.  If $r < 1-m$, then $(2)_{r,b} \Rightarrow
(1)_{r,b-r+1}$.
\item [(c)] If $r \geq -m$, then $(3)_{r,b} \Rightarrow
(4)_{r,b+r-1}$.  If $r < -m$, then $(3)_{r,b} \Rightarrow
(4)_{r,b-m}$.
\item [(d)] If $r \geq 1-m$, then $(4)_{r,b} \Rightarrow
(3)_{r,b-m}$.  If $r < 1-m$, then $(4)_{r,b} \Rightarrow
(3)_{r,b-r+1}$.
\end{itemize}
\end{lemma}

(Obviously, $(3)_{r,b} \Rightarrow (1)_{r,b}$ and $(4)_{r,b}
\Rightarrow (2)_{r,b}$, but we do not need these facts.)

\begin{proof}
As above, we write $g$ for the map $F_{s+1}X \xrightarrow{} F_{s}X$
and $g_{*}$ for the map $D_{1}^{s+1,t+1} \xrightarrow{}
D_{1}^{s,t}$, so that $D_{r}^{s,t}$ is the image of 
\[
g_{*}^{r-1} \mathcolon \pi_{t-s}F_{s+r-1}X \xrightarrow{}
\pi_{t-s}F_{s}X.
\]

(a): Assume that if $s \geq m(t-s) + b$, then
\[
g_{*}^{r-1} \mathcolon \pi_{t-s}(F_{s+r-1}X) \xrightarrow{}
\pi_{t-s}(F_{s}X)
\]
is zero; i.e., $D_{r}^{s,t} = 0$.  In the case $r \geq -m$, if $s \geq
m(t-s) + b$, then $s+r \geq m((t+r-1)-(s+r)) + b$; so we see that
$D_{r}^{s+r,t+r-1}=0$.  By the long exact sequence
\eqref{eqn-ass-exact}, we conclude that $E_{r}^{s+r-1,t+r-1} = 0$ when
$s \geq m(t-s)+b$.  Reindexing, we find that $E_{r}^{p,q} = 0$ when $p
\geq m(q-p) + b+r-1$; i.e., condition $(2)_{r,b+r-1}$ holds.  The
case $r < -m$ is similar; in this case, the long exact sequence
implies that $E_{r}^{s,t+1} = 0$.

(b): Assume that $r \geq 1-m$.  If $E_{r}^{s,t}(X) = 0$ whenever $s
\geq m(t-s) + b$, then $E_{r}^{s+r-1,t+r-2}(X) = 0$ when $s \geq
m(t-s) + b$.  So by the exact sequence \eqref{eqn-ass-exact}, we see
that $D_{r}^{s+1,t} \xrightarrow{} D_{r}^{s,t-1}$ is an isomorphism
under the same condition.  This map is induced by $g_{*} \mathcolon
\pi_{t-s-1}F_{s+1}X \xrightarrow{} \pi_{t-s-1} F_{s} X$, so we
conclude that when $s \geq m(t-s)+b$, we have
\begin{gather*}
\invlim_{q} \pi_{t-s-1} F_{q}X = D_{r}^{s,t-1}, \\
\invlim_{q}^{1} \pi_{t-s-1} F_{q}X = 0.
\end{gather*}
But by convergence of the spectral sequence, we know that $\invlim_{q}
\pi_{t-s-1} F_{q}X = 0$, so $D_{r}^{s,t-1} = \im g_{*}^{r-1} = 0$.
Reindexing gives $D_{r}^{p,q} = 0$ when $p \geq m(q+1-p)+b$; i.e.,
$(2)_{r,b}$ implies $(1)_{r,b+m}$.

If $r < 1-m$, then a similar argument shows that $D_{r}^{s-r+1,t-r+1}
= 0$.

Parts (c) and (d) are similar.
\end{proof}

It is easy to prove Theorem~\ref{thm-main}, once we have the lemma.

\begin{proof}[Proof of Theorem~\ref{thm-main}]
The proofs of the genericity of the two statements are similar, so we
only prove that condition (i) is generic.

We know by Lemma~\ref{lemma-main} that condition (i) is equivalent, up
to a reindexing, to
\begin{itemize}
\item [$(*)$] There exist numbers $r$ and $b$ so that for all $s$ and
$t$ with $s \geq m(t-s) + b$, the map $g^{r-1} \mathcolon F_{s+r-1}X
\xrightarrow{} F_{s}X$ is zero on $\pi_{t-s}$.
\end{itemize}
We show that this condition is generic.  Since the Adams tower is
functorial, if $Y$ is a retract of $X$, then the Adams tower for $Y$
is a retract of the Adams tower for $X$.  So if $F_{s+r-1}X
\xrightarrow{} F_{s}X$ is zero on $\pi_{t-s}$, then so is $F_{s+r-1}Y
\xrightarrow{} F_{s}Y$.  (Given $S^{t-s} \xrightarrow{} F_{s+r-1}Y$,
then consider
\[
\begin{CD}
S^{t-s} @>>> F_{s+r-1}Y @>>> F_{s}Y \\
@. @VV{i}V  @VV{i}V \\
@.     F_{s+r-1}X @>>> F_{s}X \\
@. @VV{j}V  @VV{j}V \\
@.     F_{s+r-1}Y @>>> F_{s}Y
\end{CD}
\]
Since $\pi_{t-s} F_{s+r-1}X \xrightarrow{} \pi_{t-s} F_{s}X$ is 0,
then the map $S^{t-s} \xrightarrow{} F_{s}X$ is null.  But $S^{t-s}
\xrightarrow{} F_{s}Y$ factors through this map, and hence is also
null.)

Given a cofibration sequence $X \xrightarrow{} Y \xrightarrow{} Z$ in
which $X$ and $Z$ satisfy conditions $(*)_{r,b}$ and $(*)_{r',b'}$,
respectively, we show that $Y$ satisfies $(*)_{r+r'-1,\max (b,
b'-r+1)}$.  Consider the following commutative diagram, in which the
rows are cofibrations:
\[
\begin{CD}
F_{s+r+r'-2} X @>>> F_{s+r+r'-2} Y @>>> F_{s+r+r'-2} Z \\
   @VVV       @VV{\alpha}V   @VV{\beta}V \\
  F_{s+r-1} X  @>>>  F_{s+r-1} Y   @>>>  F_{s+r-1} Z \\
 @VV{\gamma}V    @VV{\delta}V     @VVV \\
   F_{s}X    @>>>    F_{s}Y    @>>>   F_{s}Z
\end{CD}
\]
We assume that $s \geq m(t-s) + \max (b, b'-r+1)$, so that we have
\begin{gather*}
s \geq m(t-s) + b,\\
s + r - 1 \geq m(t-s) + b'.
\end{gather*}
If we map $S^{t-s}$ into this diagram, then since $\pi_{t-s}\beta =
0$, any map 
\[
S^{t-s} \xrightarrow{} F_{s+r+r'-2}Y \xrightarrow{\alpha}
F_{s+r-1}Y
\]
factors through $F_{s+r-1}X$.  Since $\pi_{t-s}\gamma =
0$, though, then the composite 
\[
S^{t-s} \xrightarrow{} F_{s+r+r'-2}Y \xrightarrow{\alpha} F_{s+r-1}Y
\xrightarrow{\delta} F_{s}Y
\]
is null.

This shows that condition $(*)$, and hence condition (i), is generic.
\end{proof}

The same proof, in the case $m=0$, also shows the following (using the
language of \cite{christensen-thesis}).

\begin{corollary}
If $I$ is an ideal of maps that is part of a projective class, then
the following property is generic for $E$-complete spectra
$X$\textup{:}
\begin{itemize}
\item There exist numbers $r$ and $b$ so that for all $s \geq b$, the
composite
\[
g^{r-1} \mathcolon F_{s+r-1}X \xrightarrow{} F_{s}X
\]
is in $I$.
\end{itemize}
\end{corollary}

% \bibliography{palmieri}
% \bibliographystyle{amsalpha}

\providecommand{\bysame}{\leavevmode\hbox to3em{\hrulefill}\thinspace}

\end{document}